\def\bbuildrel#1_#2^#3{\mathrel{\mathop{\kern 0pt#1}\limits_{#2}^{#3}}}

\def\NN{\mathbb N}

\def\QQ{\mathbb Q}

\def\PP{\mathbb P}

\def\RR{\mathbb R}

\def\ms{\medskip}
\def\ss{\smallskip}

\def\w{\thinspace\hbox{\hsize 14pt \rightarrowfill}\thinspace}

\def\0{\hbox{$\emptyset$}}
\def\A{\hbox{$\mathcal A$}}

\def\sub{\subseteq}

\def\Ra{\Rightarrow}
\def\Lra{\Leftrightarrow}

\documentclass[a4paper,12pt]{amsart}

\usepackage{amssymb,enumerate,mathrsfs}
\usepackage[utf8]{inputenc} 

\theoremstyle{plain}

\newcommand{\co}{\mathfrak c}

\newcommand{\sk}{[\NN]^{<\NN}}
\newcommand{\seq}{\NN^{<\NN}}

\newcommand{\cantor}{2^{\NN}}
\newcommand{\baire}{\NN^{\NN}}
\newtheorem{theorem}{Theorem}[section]

\newtheorem{lemma}[theorem]{Lemma}

\newtheorem{example}[theorem]{Example}

\newtheorem{proposition}[theorem]{Proposition}

\newtheorem{remark}[theorem]{Remark}
\newtheorem{cor}[theorem]{Corollary}

\numberwithin{equation}{section}

\begin{document}

\title{On two consequences of CH established by Sierpi\' nski}

\author{R. Pol and P. Zakrzewski}
\address{Institute of Mathematics, University of Warsaw, ul. Banacha 2,
02-097 Warsaw, Poland}
\email{pol@mimuw.edu.pl, piotrzak@mimuw.edu.pl}




\begin{abstract}

We study the relations between two consequences of the Continuum Hypothesis discovered by Wacław Sierpiński,
concerning uniform continuity of continuous functions and uniform convergence of sequences of real-valued functions, defined on subsets of the real line of cardinality continuum.

\end{abstract}

\maketitle
\section{Introduction}\label{sec:1} 
In his classical treaty \textit{Hypoth\`{e}se du continu} \cite{Si} Wacław Sierpiński distinguished the following consequences of the Continuum Hypothesis (CH) (the notation is taken from \cite{Si}):

\begin{enumerate}
	\item[${C}_8$] There exists a continuous function $f:E\to \RR$, $E\sub\RR$, $|E|=\co$, not uniformly continuous on any uncountable subset of $E$.
	
	\item[$C_9$] There is a sequence of functions $f_n:E\to\RR$, $E\sub\RR$, $|E|=\co$, converging pointwise but not converging uniformly on any uncountable subset of $E$.
\end{enumerate}

Sierpiński established the equivalences of $C_9$ to several other statements, notably, to the existence of a matrix of sets of real numbers (called in \cite{ba-ha} a {\sl BK-matrix}), constructed under CH by Banach and Kuratowski \cite{ba-ku} (statement $C_{11}$ in \cite{Si}).

Bartoszyński and Halbeisen \cite{ba-ha} (see also \cite{ha}) proved that the existence of a BK-matrix is independent of CH. They also pointed out that the existence of a BK-matrix (hence statement $C_9$) is equivalent to the existence of a subset of $\baire$ of cardinality $\co$, intersecting each compact set in $\baire$ in an at most countable set (following \cite{ba-ha} we shall call such sets {\sl K-Lusin}), see \cite[Proposition 1.1 and Lemma 2.3]{ba-ha}, cf. also \cite{ro}.

Sierpiński \cite{si-2} noticed that $C_8$ implies ${C}_9$ but he did not discuss the converse implication. However,
in \textit{Topology I} by Kuratowski \cite{kur}, footnote (3) on page 533 suggests that the two  statements are in  fact equivalent. We are not aware of any publication addressing the implication $C_9 \Ra C_8$ and this note is the result of our pondering on this matter.


\ms 

We shall consider the following stratifications of statements $C_8$ and $C_9$ for  uncountable cardinals  $\kappa\leq\lambda\leq\co$: 
\begin{enumerate}
	\item[$C_8(\lambda,\kappa)$] 
	There exist a set $E\sub \RR$ 
	of cardinality $\lambda$ and a continuous function $f:E\w \RR$, which is not uniformly continuous on any subset of $E$ of cardinality $\kappa$.
	
		\item[$C_9(\lambda,\kappa)$] There exists a set $E\sub \RR$ of cardinality $\lambda$ (equivalently: for any set $E\sub \RR$ of cardinality $\lambda$)
	and	there is a sequence of functions $f_n:E\to\RR$,  converging on $E$ pointwise but not converging uniformly on any subset of $E$ of cardinality $\kappa$.
\end{enumerate}

Clearly, statements $C_i$  are $C_i(\co,\aleph_1)$ in our notation, and $C_i$ implies $C_i(\lambda,\kappa)$ for all uncountable cardinals  $\kappa\leq\lambda\leq\co$, $i=8,9$.

\ms

In this note we prove (in ZFC) that, in particular:

\begin{itemize}
	
\item $C_8(\co,\co)\Lra C_9(\co,\co)$ provided that the cardinal $\co$ is regular, and under this assumption each of these statements is equivalent to the assertion $\mathfrak{d}=\co$ (cf. Corollary \ref{C_8_iff_C_9}),

\item $C_8(\aleph_1,\aleph_1)\Lra C_9(\aleph_1,\aleph_1)$,
and each of these statements is equivalent to the assertion $\mathfrak{b}=\aleph_1$ (cf. Theorem \ref{C_8(w_1,w_1)_iff_C_9(w_1,w_1)}).
\end{itemize}
Here $\mathfrak{d}$ and $\mathfrak{b}$ denote, as usual, the smallest cardinality of a dominating and, respectively, an unbounded family in $\baire$ corresponding to the ordering of eventual domination $\leq^*$ (cf. \cite{ha}).  

\ms 

An important role in our considerations is played by the notion of a K-Lusin set which we  extend  (cf.  \cite{ba-ha}) declaring that an uncountable subset $E$ of a Polish space $X$ is a {\sl $\kappa$-$K$-Lusin set in $X$},  $\aleph_1\leq \kappa\leq\co$, if $|E\cap K|<\kappa$  for every  compact set $K\sub X$.

\ms 

 The existence of a $\kappa$-$K$-Lusin set of cardinality $\lambda$ in $\baire$ is equivalent to $C_9(\lambda,\kappa)$ (cf. Theorem \ref{C_9}) and if $E\sub \RR$, $|E|=\lambda$, is a witnessing set for $C_8(\lambda,\kappa)$, then $E$ is a $\kappa$-$K$-Lusin set in some $G_\delta$-extension of $E$ (cf. Proposition \ref{C_8->C_9}).
 
  However, it is not the case that every $\kappa$-$K$-Lusin set is a witnessing set for $C_8(\lambda,\kappa)$. In particular, assuming CH, we show that there is a $K$-Lusin set in the irrationals of cardinality $\co$ such that  every continuous function $f:E\w \RR$ is uniformly continuous on an uncountable subset of $E$ (cf. Theorem \ref{K-Lusin}). Our reasoning to that effect yields also that  for every continuous function $f:X\to\RR$ defined on a $G_\delta$-set $X$ in the irrationals, there exists a closed copy of irrationals $P$ in $X$ such that $f$ is uniformly continuous on $P$  (cf. Theorem \ref{unif_cont_on_P}).


\ms 
The paper is organized as follows.
\smallskip

In Section \ref{sec:2} we establish the equivalences of $C_9(\lambda,\kappa)$ to several other statements, notably, to its topological counterparts (see Theorem \ref{C_9}).

\smallskip 
Section \ref{sec:3} is devoted to $C_8(\lambda,\kappa)$ and its relations to  $C_9(\lambda,\kappa)$ including proofs of the equivalence $C_8(\lambda,\kappa)\Lra C_9(\lambda,\kappa)$ for $\kappa=\lambda=\co$ (if the cardinal $\co$ is regular) and $\kappa=\lambda=\aleph_1$. 
We end this section by listing some additional set-theoretic assumptions under which the equivalence $C_8(\co,\aleph_1)\Lra C_9(\co,\aleph_1)$ is true. 
Although the status of the implication $C_9(\co,\aleph_1)\Ra C_8(\co,\aleph_1)$ remains unclear, these observations point out at difficulties in refuting it.

\smallskip

In Section \ref{sec:4} we gathered some comments and additional results related to the topic without proofs -- we plan to present details  elsewhere.

\ms

In this note $\PP$ always denotes the set of irrationals of the unit interval $[0,1]$. It is homeomorphic to the Baire space $\baire$, the countable product of the set of natural numbers $\NN=\{0,1,2,\ldots\}$ with the discrete topology (cf. \cite{ke}).

\section{Uniform convergence of pointwise convergent sequences of functions and statement $C_9(\lambda,\kappa)$}\label{sec:2} 

The following result is based on Sierpiński's reasoning \cite{si-1}, cf. Remark \ref{Sierp}(1)
(an extension of this result is formulated in Section  \ref{sec:4.3}).

\begin{theorem}\label{seq_on_Polish}
For any Polish space $X$ there is a sequence $f_1\geq f_2\ldots$ of continuous functions $f_n:X\w [0,1]$ which converges to zero pointwise but does not converge uniformly on any set with non-compact closure in $X$. 
\end{theorem}

\begin{proof}

Since $X$ embeds as a closed subspace in $[1,+\infty)^{\NN_+}$ (cf. \cite[Theorem 4.17]{ke}), 
where $\NN_+=\{n\in\NN:\ n>0\}$,
 it is enough to construct desired functions on $[1,+\infty)^{\NN_+}$. So, with no loss of generality, we simply assume that $X=[1,+\infty)^{\NN_+}$. 
 
 We begin with the Sierpiński functions  $s_n:X\w \{0,1,\frac{1}{2},\frac{1}{3},\ldots\}$, $n=1,2,\ldots $, defined by (cf. Remark \ref{Sierp})
	\begin{enumerate}
	\item[(1)] $s_n(x)=\begin{cases} \frac{1}{\min\{i:\ x(i)\geq n\}}\quad&\hbox{if}\quad  x(\NN_+)\cap [n,+\infty)\neq\emptyset,\\
	0 &\hbox{if} \quad x(\NN_+)\sub [0,n).
\end{cases}$
\end{enumerate}

We shall check that 
\begin{enumerate}
	\item[(2)] $s_1\geq s_2\geq\ldots$ and $\lim\limits_{n\to\infty}s_n(x)=0$ for every $x\in X$,
	\item[(3)] for any $A\sub X$, if the sequence $(s_n)_{n=1}^{\infty}$ converges uniformly on $A$, then the closure $\bar{A}$ is compact. 
\end{enumerate}

The monotonicity in (2) is clear. If $x\in X$ and $p\in \NN_+$ is given, then for any $n>\max\{x(i):\ i\leq p\}$ we have $s_n(x)<\frac{1}{p}$, which gives the second part of (2).

\smallskip 

To make sure that (3) is true, we shall follow closely Sierpiński \cite{si-1}. Let $A\sub X$ and assume that the sequence $(s_n)_{n=1}^{\infty}$ converges uniformly on $A$. This means that for each $i\in \NN_+$ there is $\varphi(i)\in \NN_+$ such that 
$s_m(x)<\frac{1}{i}$, whenever $m\geq \varphi(i)$ and $x\in A$.

By (1), for any $x\in X$ and $i\in\NN_+$ we have $s_{\lfloor x(i)\rfloor}(x)\geq \frac{1}{i}$, and hence $x(i)\leq \varphi(i)$, for any $x\in A$. Therefore, $A$ is contained in the compact set $\prod\limits_{i=1}^\infty [1,\varphi(i) ]\sub X$, and hence its closure $\bar{A}$ in $X$ is compact.

\smallskip

Let us verify that for each $n\in\NN_+$
\begin{enumerate}
	\item[(4)] the function $s_n$ is upper-semicontinuous,
\end{enumerate}
i.e., for any $r>0$ the set $\{x\in X:\ s_n(x)<r\}$ is open in $X$. Since $s_n$ is bounded by 1, it is enough to consider $r\leq 1$.

So let us fix $n\in\NN_+$, $r\leq 1$ and $a\in X$ with $s_n(a)<r$, and for any $p\in\NN_+$, let us consider the open set $V_p$ defined by
 \begin{enumerate}
 	\item[(5)] $V_p=\{x\in X:\ x(i)<n\quad\hbox{for all}\quad i\leq p \}$.
 \end{enumerate}
We shall show that we can always find $p$ such that $V_p$ is a neighbourhood of $a$ contained in the set $\{x\in X:\ s_n(x)<r\}$.

\smallskip 

If $s_n(a)=0$, i.e., $a(i)<n$ for all $i\in \NN_+$, cf. (1), then taking  $p$ such that $\frac{1}{p}<r$, we have $a\in V_p$ and $s_n(x) <\frac{1}{p}<r$ for every $x\in V_p$, cf. (1) and (5).

If $s_n(a)=\frac{1}{m}$, where $m=\min\{i:\ a(i)\geq n\}$, then since $s_n(a)<1$, we have $a(1)<n$. Hence  $m\geq 2$ and let $p=m-1$. Then $p\geq 1$, $a\in V_p$ and for any $x\in V_p$, $s_n(x)\leq \frac{1}{m}<r$.

\smallskip 

Having checked (4), we apply a classical theorem of Hahn (cf. \cite[1.7.15(c)]{eng}) to get, for each $n$, continuous functions $f_{n,i}:X\w [0,1]$, $i=1,2,\ldots$, such that 
\begin{enumerate}
	\item[(6)] $f_{n,1}\geq f_{n,2}\geq\ldots$ and $\lim\limits_{i\to\infty} f_{n,i}(x)=s_n(x)$ for every $x\in X$.
\end{enumerate}

Finally, we define
\begin{enumerate}
	\item[(7)] $f_n(x)=\min\limits_{i,j\leq n} f_{i,j}(x)$ for $x\in X$.
\end{enumerate}

Clearly, the sequence $f_1\geq f_2\geq\ldots$ consists of continuous functions and converges pointwise to zero. Moreover, $f_n(x)\geq s_n(x)$ for any $n\in\NN_+$ and $x\in X$. Consequently, for any $A\sub X$, if the sequence $(f_n)_{n=1}^{\infty}$ converges uniformly on $A$, then so does the sequence $(s_n)_{n=1}^{\infty}$ and hence by (3), $\bar{A}$ is compact.

	\end{proof}

\begin{remark}\label{Sierp} \hfill\null
	\begin{enumerate}
		
	\item  The original Sierpiński functions were defined on $\NN_+^{\NN_+}$ by the formula:
	$$s_n(x)=\begin{cases} \frac{1}{\min\{i:\ x(i)= n\}}\quad&\hbox{if}\quad  n\in x(\NN_+),\\
		0 &\hbox{if} \quad n\notin x(\NN_+).
	\end{cases}$$
		Sierpiński was interested  in neither regularity of the functions (in fact, $s_n$ are continuous on $\NN_+^{\NN_+}$) nor  the monotonicity of the function sequence.
	
	\smallskip
	\item An approach similar to Sierpiński's idea, in a different setting, was rediscovered by  Pincirolli \cite[Lemma 2 and Proposition 7]{pinc}. 
		\end{enumerate}
	\end{remark}

\smallskip
Recall that $C_9(\lambda,\kappa)$ abbreviates the following statement:
\begin{enumerate}

	\item[] There exists a set $E\sub \RR$ of cardinality $\lambda$ (equivalently: for any set $E\sub \RR$ of cardinality $\lambda$) and
there is a sequence of functions $f_n:E\to\RR$,  converging on $E$ pointwise but not converging uniformly on any subset of $E$ of cardinality $\kappa$.
	
\end{enumerate}

The following result provides some topological counterparts to $C_9(\lambda,\kappa)$.

\begin{theorem}\label{C_9}
For any uncountable cardinals $\kappa\leq\lambda\leq\co$	the following are equivalent:
	\begin{enumerate}
		\item $C_9(\lambda,\kappa)$,
		 \smallskip
		 
		\item there is a set $A\sub \baire$ of cardinality $\lambda$ and a sequence $g_1\geq g_2\ldots$ of continuous functions  $g_n:A\to \RR$,
		  which converges to zero pointwise but does not converge uniformly on any set of cardinality $\kappa$,
		  
		\smallskip 
		
		\item 
		there is a $\kappa$-K-Lusin set of cardinality $\lambda$ in $\baire$,
		
		\smallskip 
		\item 
		there is a Polish space $X$ and a $\kappa$-K-Lusin set of cardinality $\lambda$ in $X$.
		 
		\end{enumerate}
\end{theorem}

\begin{proof}\hfill\null
	

	
		$(1)\Ra (2)$. Subtracting from each function in $C_9(\lambda,\kappa)$ the limit function, we get a sequence $f_n:E\w \RR$, $|E|=\lambda$,
			which converges to zero pointwise but does not converge uniformly on any set of cardinality $\kappa$. For every $n\in\NN_+$ and $x\in E$ let $u_n(x)=\max\{|f_i(x)|:\ i\geq n\}$ (recall that $\lim\limits_{i\to\infty}f_i(x)=0$, hence the maximum is attained). Let us note that $u_1\geq u_2\ldots$ and  $0\leq |f_n|\leq u_n$ for each $n$. The properties of the sequence $(f_n)_{n=1}^\infty$ yield readily that 	 the sequence $(u_n)_{n=1}^\infty$ converges to zero pointwise on $E$ but it does not converge uniformly on any subset of $E$  of cardinality $\kappa$.

		 Let $h:E\w A$ be a bijection onto a set $A\sub \cantor$ such that all the functions $g_n=u_n\circ h^{-1}$ are continuous (one may define $h$ as the Marczewski characteristic function (cf. \cite{mi}) of a countable family $\{E_n:\ n\in\NN\}$ of subsets of $E$, separating the points of $E$ and containing all sets of the form $u_n^{-1}\bigl((p,q)\bigr)$, where $n\in\NN_+$ and $p<q$ are rationals). Then the sequence $g_1\geq g_2\ldots$ of continuous functions  $g_n:A\to \RR$ is as required.
	
	\smallskip
	
	$(2)\Ra (3)$. Let us fix a set $A\sub \baire$ of cardinality $\lambda$ and a sequence $g_1\geq g_2\ldots$ of continuous functions  $g_n:A\w \RR$,
	which converges to zero pointwise but does not converge uniformly on any subset of $A$ of cardinality $\kappa$. Let $H$ be a $G_\delta$-set in $\baire$ with $A\sub H\sub \bar{A}$ and such that each $g_n$ extends to a continuous function $\tilde{g}_n: H\w \RR$. Since $\tilde{g}_1\geq \tilde{g}_2\geq\ldots$, for any $x\in H$ we have $$
	\lim\limits_{n\to\infty}\tilde{g}_n(x)=0\ \Lra\ 
	\forall p\in\NN_+\exists n\in\NN_+\ \tilde{g}_n(x)<\frac{1}{p},
	$$
	so  the set $G=\{x\in H:\ \lim\limits_{n\to\infty}\tilde{g}_n(x)=0\}$ is a $G_\delta$-set in $\baire$ containing $A$. Now, if $K\sub G$ is compact, then by the Dini theorem (see \cite[Lemma 3.2.18]{eng}), the sequence $(\tilde{g}_n)_{n=1}^\infty$ converges uniformly on $K$, hence also $({g}_n)_{n=1}^\infty$ converges uniformly on $A\cap K$, and therefore $|A\cap K|<\kappa$.
	
	Finally, let $w:G\w \baire$ embed $G$ onto a closed subspace of $\baire$ (see \cite[Theorem 7.8]{ke}) and let $E=w(A)$. Then $|E\cap K|<\kappa$  for every  compact
	set $K\sub\baire$, as required.
	
	\smallskip 
	
	$(3)\Ra (4)$. This implication is trivial.
	
	\smallskip 
	
	$(4)\Ra (1)$.  Let  $E$ be a subset of cardinality $\lambda$ of a Polish space $X$ such that $|E\cap K|<\kappa$  for every  compact
	set $K\sub X$. Theorem \ref{seq_on_Polish} provides us with a sequence $f_1\geq f_2\ldots$ of continuous functions $f_n:X\w [0,1]$ which converges to zero pointwise but does not converge uniformly on any set with non-compact closure in $X$. Consequently, the sequence $({f}_n)_{n=1}^\infty$ converges to zero pointwise on $E$ but any set $M\sub E$ of cardinality $\kappa$ has a non-compact closure in $X$, so $({f}_n)_{n=1}^\infty$ does not converge uniformly on $M$. Clearly, this completes the proof of (1).
	
	\end{proof}

In two important cases statements $C_9(\kappa,\lambda)$ are characterized in terms of basic cardinal characteristics of the continuum, cf. \cite{bla}. 

\begin{cor}\label{b,d}\hfill\null
	
		\begin{enumerate}
		\item $C_9(\aleph_1,\aleph_1)\Lra  \mathfrak{b}=\aleph_1$,
 \item $C_9(\co,\co)\Lra  \mathfrak{d}=\co$, provided that the cardinal $\co$ is regular (more precisely, $\mathfrak{d}=\co$ implies $C_9(\co,\co)$ under no additional assumptions on $\co$).
	\end{enumerate}
\end{cor}

\begin{proof} We shall repeatedly make use of Theorem \ref{C_9}.
	
	\smallskip 
	
(1).  Assume $C_9(\aleph_1,\aleph_1)$ and let  $E\sub \baire$ be a set of cardinality $\aleph_1$  whose intersection with every  compact
set $K\sub\baire$ is countable. Clearly, $E$ is unbounded in $(\baire,\leq^*)$, hence $\mathfrak{b}=\aleph_1$.

Conversely, if $\mathfrak{b}=\aleph_1$, then  any subset of 
$\baire$ of the form $\{f_\alpha:\alpha< \mathfrak{b}\}$, 
where
\begin{itemize}
	\item $\alpha<\beta<\mathfrak{b}$ implies $f_\alpha<^*f_\beta$,
	\item for every $f\in\baire$ there is $\alpha<\mathfrak{b}$ with $f_\alpha\nleq^*f$,
\end{itemize}
has countable intersection with every  compact $K\sub\baire$. 

\smallskip 

(2). Assume $C_9(\co,\co)$ and let  $E\sub \baire$ be a set of cardinality $\co$ such that $|E\cap K|<\co$ for  every  compact
set $K\sub\baire$. Let $\{g_\alpha:\alpha< \mathfrak{d}\}$ be a dominating set in $\baire$. In particular $E=\bigcup_{\alpha<\mathfrak{d}} \{f\in E:\ f<^* g_\alpha\}$ and the regularity of $\co$ implies that $\mathfrak{d}=\co$.

Conversely, if $\mathfrak{d}=\co$, then
 any subset of 
$\baire$ of the form $\{g_\alpha:\alpha< \co\}$, 
where
\begin{itemize}
	\item $\alpha<\beta<\co$ implies $g_\beta\not\leq^*g_\alpha$,
	\item for every $f\in\baire$ there is $\alpha<\co$ with $f <^*g_\alpha$,
\end{itemize}
has the property that $|E\cap K|<\co$ for every  compact $K\sub\baire$.

\end{proof}

The proof of $(4)\Ra (1)$ in Theorem \ref{C_9} yields also the following result.

\begin{cor} For any uncountable cardinals $\kappa\leq\lambda\leq\co$,
if	$E$ is a $\kappa$-K-Lusin set of cardinality $\lambda$ in   a Polish space $X$,  
 then there exists a sequence $f_1\geq f_2\ldots$ of continuous functions  $f_n:E\to \RR$,
which converges to zero pointwise but does not converge uniformly on any subset of $E$ of cardinality $\kappa$.
\end{cor}


\begin{remark}
As was mentioned in the introduction, the notion of a $K$-Lusin set was introduced by Bartoszyński and Halbeisen \cite{ba-ha}, where it was pointed out that a reasoning of Banach and Kuratowski \cite{ba-ku}, establishing under CH the existence of  a BK-Matrix, actually shows that  the existence of a BK-matrix is equivalent to the existence of a $K$-Lusin set of cardinality $\co$. Earlier, Sierpiński \cite{si-1} proved that a BK-Matrix exists if and only if $C_9$ holds. Combining these two results, we get the equivalence "$C_9 \Lra\ \hbox{there exists a K-Lusin set of cardinality $\co$}$`` which was obtained in Theorem \ref{C_9} by a different reasoning.  

	\end{remark}

\section{Uniform continuity of continuous functions and statement $C_8(\lambda,\kappa)$}\label{sec:3} 

Recall that $C_8(\lambda,\kappa)$ stands for the following statement:
\begin{enumerate}
	\item[] 
There exists a set $E\sub \RR$ 
  of cardinality $\lambda$ and a continuous function $f:E\w \RR$, which is not uniformly continuous on any subset of $E$ of cardinality $\kappa$.
\end{enumerate}


Sierpiński \cite{si-2} proved that $C_8$ implies $C_9$ and his argument can be easily adapted to establish a more general implication concerning $C_i(\lambda,\kappa)$. Instead of repeating the argument of Sierpiński we present a proof based on Theorem \ref{C_9} which gives some additional information about the involved sets.

\begin{proposition}\label{C_8->C_9}
For any uncountable cardinals  $\kappa\leq\lambda\leq\co$:
$$
C_8(\lambda,\kappa)\Ra C_9(\lambda,\kappa).
$$
Moreover, if a set $E\sub \RR$, $|E|=\lambda$, together with a continuous function $f:E\w \RR$ witness $C_8(\lambda,\kappa)$, then there is a $G_\delta$-set $G$ in $\RR$ such that $E\sub G$ and 
$E$ is a $\kappa$-$K$-Lusin set in $G$.
\end{proposition}

\begin{proof}
Let us extend $f$ to a continuous function $\tilde{f}:G\to\RR$ over a $G_\delta$-set $G\sub\RR$. 
Now, if $K\sub G$ is compact, then the extension $\tilde{f}$ is uniformly continuous on $K$, hence so is $f$ on $E\cap K$. Therefore, $|E\cap K|<\kappa$, as
$f$ is not uniformly continuous on any set of cardinality $\kappa$. This shows that the equivalent to $C_9(\lambda,\kappa)$ statement, formulated in  Theorem \ref{C_9}(4), is true, completing the proof.
\end{proof}

In the rest of this note we investigate the possibility of reversing the above implication, at least for some pairs of uncountable cardinals $\kappa\leq\lambda\leq\co$. 

\smallskip 

In view of Proposition \ref{C_8->C_9}, a related question is whether a $\kappa$-$K$-Lusin set $E$  in $\PP$
 always carries a continuous function $f:E\w \RR$, which is not uniformly continuous on any set of cardinality $\kappa$. The negative answer (cf. Theorem \ref{K-Lusin}) is a consequence of the following general result,
closely related to the ``limit systems'' of Hurewicz \cite{hu}.

\begin{theorem}\label{unif_cont_on_P}
Let $X$ be a Polish non $\sigma$-compact space and let $d$ be a compatible completely bounded metric on $X$. Then for every continuous function $f:X\to\RR$ there exists a closed copy of irrationals $P$ in $X$ such that $f$ is uniformly continuous on $P$ in the metric $d$. 	
	
\end{theorem}

\begin{proof}
Let $(\hat{X},\hat{d})$ be the completion of $(X,d)$; then $\hat{X}$ is compact, $d$ being totally bounded. Since $X$ is not $\sigma$-compact, by a theorem of Hurewicz (see \cite[Theorem 7.10]{ke}), $X$ contains a closed in $X$ copy of the irrationals $G$. Let $\rho$ be a complete metric on $G$. 

We shall use generalized Hurewicz systems in the setting considered in  \cite[Section 2.4]{p-z-1} and \cite[Section 2]{p-z-2}. Namely, we shall define  a pair of families: 
$(U_s)_{s\in\seq}$ of subsets of $G$, and $(x_s)_{s\in\seq}$ of points in $\hat{X}$   with the following properties (the closures are taken in $\hat{X}$, $B_{\hat{d}}(x_s,\varepsilon)=\{x\in\hat{X}:\hat{d}(x_s,x)<\varepsilon\}$ and for $A\sub G$, $\hbox{\rm diam}_{\rho} (A)$ or $\hbox{\rm diam}_{d} (A)$ stand for the diameter with respect to $\rho$ or $d$):

  \begin{enumerate}
	\item $U_s$ is relatively open in $G$, $U_s\neq\emptyset$, 
	
	\smallskip 
	
	\item  $\hbox{\rm diam}_{\rho} (U_s) \leq 2^{-length(s)}$,
	
	\smallskip

	
	
	\item $\overline{U_s}\cap \overline{U_t}=\emptyset$ for distinct $s,\ t$ of the same length,
	
	\smallskip 
	
	\item $\overline{U_{s \hat\ i}}\cap G\sub U_s$,
	
	\smallskip 
	
	

	\item $x_s \in \overline{U_s}\setminus G$,
	
	\smallskip
	
	\item $x_s\not\in \overline{U_{s  \hat\ i}}$ for any $i\in\NN$, 
	
	\smallskip 
	
	\item  each neighbourhood of $x_s$ in $\hat{X}$ contains all but finitely many $U_{s\hat\ i}$; in particular, $\lim_i\hbox{\rm diam}_d (U_{s\hat\ i})=0$,
	
	\smallskip
	
		\item $\hbox{\rm diam}(f(U_{s\hat\ i}))\leq 2^{i}$ for any $i\in\NN$,
	
	
	\smallskip
	
	\item if $c_i \in U_{s\hat\ i}$ for each $i\in\NN$, then the sequence $(f(c_i))_{i\in\NN}$ is convergent in $\RR$.

\end{enumerate} 

To define sets $U_s$ and points $x_s$ we proceed as follows.

\smallskip 

 Let $U_{\emptyset}$ be a non-empty relatively open set in $G$ such that $f$ is bounded on $U_{\emptyset}$ and $\hbox{\rm diam}_{\rho} (U_{\emptyset}) \leq 1$.

\smallskip 

At the inductive step let $n\geq 0$ and assume that we have already defined $U_s$ and $x_t$ for $s\in [\NN]^{\leq n}$ and $t\in[\NN]^{< n}$ satisfying the required conditions. Fix $s$ with $length(s)=n$ and  pick $x_s\in \overline{U_s}\setminus G$  arbitrarily (this is possible since $G$ does not contain compact sets with non-empty interior).  Let us choose points $a_n\in U_s$ such that $\lim_n a_n = x_s$ and the sequence $(f(a_n))_n$ is convergent (first, we choose $b_n\in U_s$ so that $\lim_n b_n = x_s$ and then, using the fact that the sequence $(f(b_n))_n$ is bounded, we choose its convergent subsequence). Next, using the continuity of $f$ on $G$, let us enlarge each $a_n$ to its open neighbourhood $U_{s\hat\  n}$ in $G$ so that relevant instances of conditions (1)--(8)  are satisfied. Then (8) and the fact that the sequence $(f(a_n))_n$ is convergent readily yield (9). 

Let
\begin{enumerate}
	\item[(10)] $P = \bigcap_n \bigcup \{U_s: \hbox{length(s)} = n\}.$
\end{enumerate}
be the copy of the irrationals  determined by the generalized Hurewicz system  $(U_s)_{s\in\seq}$, $(L_s)_{s\in\sk}$ (\cite[Section 2]{p-z-2}), where $L_s=\{x_s\}$ for each $s\in \seq$. 
In particular,   $\overline{P} = P \cup \{x_s: s\in\seq\}$, so  $P=\overline{P}\cap G$ is closed in $G$, and hence also in $X$.

\smallskip 

We claim that for each $s\in\seq$

\begin{enumerate}
	\item[(11)] $\inf\limits_{\varepsilon > 0} \hbox{\rm diam} \Bigl(f\bigl(B_{\hat{d}}(x_s,\varepsilon)\cap P\bigr)\Bigr) =0.$
\end{enumerate}

To justify the claim, let us fix $s\in\seq$
and  for each  $i\in\NN$ let us pick $c_i\in U_{s\hat\ i}$. By (9), $\lim\limits_{i\to\infty} f(c_i)=r$, and let $J$ be an arbitrary open interval containing $r$. From (7) and (8) we get $i_0$ such that $f(U_{s\hat\ i})\sub J$, whenever $i>i_0$.

Now, 
we can find an 
$\varepsilon > 0$ such that $B_{\hat{d}}(x_s,\varepsilon)$ is disjoint from any $U_t$ with $t\neq s$ and $length(t)=length(s)$, i.e., cf. (10),
\begin{enumerate}
	\item[(12)]
	$W= B_{\hat{d}}(x_s,\varepsilon)\cap P=B_{\hat{d}}(x_s,\varepsilon)\cap U_s$.

\end{enumerate}

For suppose that for every $\varepsilon>0$, (12) is false. This allows us to define a sequence $(z)_{n\in\NN}$ converging to $x_s$  such that the set $Z=\{z_n:n\in\NN\}$ is disjoint from $U_s$.

By (3) and (5), $Z\cap U_t$ is finite for any $t\in\seq$ with $t\neq s$ and $length(t)=length(s)$.  Then, since $Z\sub U_\emptyset$, it follows that  we can find $t\in\seq$ such that $length(t)<length(s)$, $Z\cap U_t$ is infinite but $Z\cap U_{t\hat\ i}$ is finite for each $i\in\NN$.  By (9), $\lim\limits_{n\to \infty} z_n=x_t$, however, by (3) and (6), $x_t\neq x_s$, and this contradiction completes the justification of (12).

Next, by appealing to (6), we can make $\varepsilon$ still smaller to ensure that $B_{\hat{d}}(x_s,\varepsilon)$ omits also all $U_{s\hat\ i}$ with $i\leq i_0$.
Consequently, cf. (10) and (12), $W\sub\bigcup\limits_{i>i_0}U_{s\hat\ i}$, hence $f(W)\sub J$. 

\smallskip 

This completes the justification of the claim, and let us note that (11) means exactly that the oscillation of $f$ at any point of $\overline{P}\setminus P=\{x_s: s\in\seq\}$ is zero. This guarantees that $f$ can be extended continuously over $\overline{P}$. By compactness of $\overline{P}$ this extension is uniformly continuous, and in effect  we get uniform continuity of $f$ on $P$. 

\end{proof}

As a corollary we obtain  that not every $\kappa$-$K$-Lusin set of cardinality $\lambda$ in the irrationals is a witnessing set for $C_8(\lambda,\kappa)$ (cf. \cite[a remark at the end of the paper]{si-2}).

\begin{theorem}\label{K-Lusin}
If $\mathfrak{d}=\co$, then there is a $\co$-$K$-Lusin set $E$ in $\PP$ of cardinality $\co$
 such that every continuous function $f:E\w \RR$ is uniformly continuous on a subset of $E$ of cardinality $\co$.	

In particular, assuming CH, there is a $K$-Lusin set $E\sub\PP$ of cardinality $\co$ such that  every continuous function $f:E\w \RR$ is uniformly continuous on an uncountable subset of $E$.		
\end{theorem}

\begin{proof}
	 We list all compact sets in $\PP$ as $( K_\alpha:\ \alpha<\co)$, and all closed copies of irrationals in $\PP$ as  $( P_\alpha: \ \alpha<\co)$, where each closed copy of irrationals $P$ in $\PP$ appears in this transfinite sequence $\co$-many times.
	 
	 Then we inductively pick
	 $$
	 x_\alpha\in P_\alpha\setminus\Bigl(\bigcup_{\beta<\alpha}K_\beta \cup \{x_\beta: \beta<\alpha\} \Bigr),
	 $$
	 the choice being made possible by the assumption  $\mathfrak{d}=\co$ which means  that $\PP$ is not covered by any collection of less that $\co$-many its compact subsets  (cf. \cite[Theorem 2.8]{bla}).
	 
	 We let $E=\{x_\alpha:\alpha<\co\}$. Let us notice that $E$ is a $\co$-$K$-Lusin set $E$ in $\PP$ and 
	  \begin{enumerate}
	  	\item $E$ intersects  each closed copy of  irrationals in $\PP$ in a set of cardinality $\co$. 
	  \end{enumerate}

	 To see that $E$ is a set we are looking for, let $f:E\to \RR$ be a continuous function, and let $X$ be a $G_\delta$-set in $\PP$ containing $E$ such that $f$ extends to the (uniquely defined) continuous function $\hat{f}:X\to \RR$. 
	 
	 Now, $E$ being a $\co$-$K$-Lusin set in $\PP$, it cannot be covered by countably many compact sets in $\PP$. Consequently, $X$ is a non $\sigma$-compact Polish space contained in $[0,1]$, so by Theorem \ref{unif_cont_on_P}  there exists a closed copy of irrationals $P$ in $X$ such that $\hat{f}$ is uniformly continuous on $P$ (in the metric inherited from $[0,1]$). Shrinking $P$, if necessary, we may assume that $P$ is closed also in $\PP$. By (1), $|E\cap P|=\co$ and we conclude that $f$ is uniformly continuous on 	$E\cap P$.
	 	 
\end{proof}


On the other hand, we have the following result.

\begin{theorem}\label{d=c implies C_8} 
	If $\mathfrak{d}=\co$, then there is a $\co$-$K$-Lusin  set $E$ in $\PP$ of cardinality $\co$ 
	 and a continuous function $f:E\w \RR$, which is not uniformly continuous on any set  of cardinality $\co$.
	 
	 In particular, we have $\mathfrak{d}=\co\Ra C_8(\co,\co)$.
\end{theorem}

Our proof will be based on the following proposition.

\begin{proposition}\label{covering}
Let $f:\PP\to [0,1]$ be a continuous function such that the closure $\overline{G(f)}$	in $[0,1]^2$ of the graph $G(f)$ of $f$ intersects each $\{q\}\times [0,1]$, $q\in\QQ\cap [0,1]$, in an uncountable set. Then $\PP$ cannot be covered by less than $\mathfrak{d}$ sets on which $f$ is uniformly continuous.

\end{proposition} 

\begin{proof}[Proof of Proposition \ref{covering}]
Let $\A$ be a collection of subsets of $\PP$ such that $|\A|<\mathfrak{d}$ and $f$ is uniformly continuous on each $A\in\A$. We will show that $\PP\setminus \bigcup\A\neq\emptyset$. 

Let $A\in\A $. Then $f|A$ being uniformly continuous extends continuously over the closure of $A$ in $[0,1]$, and let $K_A$ be the graph of this extension.

For each $q\in\QQ\cap [0,1]$, $K_A$ intersects $\{q\}\times [0,1]$ in at most a singleton, and hence   $|\bigcup\limits_{A\in {\mathcal A}}\bigl(K_A\cap(\{q\}\times [0,1])\bigr)|<\mathfrak{d}$. If $V$ is an open neighbourhood of $(t,f(t))$, $t\in\PP$, in $[0,1]^2$, then there are non-empty open intervals $I_1, \ I_2$ in $[0,1]$ such that $t\in I_1$, $I_1\times I_2\sub V$  and $f(I_1\cap\PP)\sub I_2$. It follows that for every $q\in I_1\cap\QQ$, 
$\overline{G(f)}\cap (\{q\}\times \bar{I}_2)=\overline{G(f)}\cap (\{q\}\times [0,1])$, so by the properties of $f$,
$|\overline{G(f)}\cap (\{q\}\times I_2)|=\co$. 
Consequently, the set
$$
H=\overline{G(f)}\cap (\QQ\times [0,1])\setminus \bigcup\limits_{A\in {\mathcal A}}K_A
$$
is dense in $\overline{G(f)}$.

Since $G(f)$ is a $G_\delta$-set dense in $\overline{G(f)}$, by the Baire theorem, each $F_\sigma$-set covering $G(f)$ must hit $H$, and by the Kechris-Louveau-Woodin theorem (see \cite[Theorem 21.22]{ke}), we get a Cantor set $C\sub G(f)\cup H$ such that $P=C\cap G(f)$ is a copy of the irrationals, closed in $G(f)$.

For each $A\in\A $, $K_A$ being compact, the set $K_A\cap P=K_A\cap C$ is compact and since $|\A|<\mathfrak{d}$ it follows that $P\setminus \bigcup\limits_{A\in {\mathcal A}}K_A\neq \emptyset$  (cf. \cite[Theorem 2.8]{bla}). In effect, $G(f) \setminus \bigcup\limits_{A\in {\mathcal A}}K_A\neq \emptyset$ which proves that $\PP\not\subseteq \bigcup \A$.
\end{proof}


\smallskip 

With Proposition \ref{covering} in hand, we can easily get 
Theorem \ref{d=c implies C_8}.

\begin{proof}[Proof of Theorem \ref{d=c implies C_8}]
	Let $f:\PP\to [0,1]$ be as in Proposition \ref{covering} 
	 (see Example \ref{kur-sie}).
	
	Let us list as $(F_\alpha:\ \alpha<\co)$ all closed sets in $\PP$ on which $f$ is uniformly continuous. Since $\mathfrak{d}=\co$, by the assertion of Proposition \ref{covering}, we can inductively pick points 
	$$
	x_\alpha\in \PP\setminus\Bigl(\bigcup\limits_{\beta<\alpha} F_\beta \cup \{x_\beta:\beta<\alpha \}\Bigr), \alpha<\co,
	$$
	and finally let $E=\{x_\alpha:\alpha<\co\}$.

	Then if $A\sub E$ and $f$ is uniformly continuous on $A$, the closure of $A$ in $\PP$ is listed as some $F_\alpha$, and hence $|A|<\co$.
	
	Likewise, every compact set $K$ in $\PP$ is on the list, so $|E\cap K|<\co$ which shows that $E$ is a $\co$-$K$-Lusin set in $\PP$.
	
\end{proof}

For the sake of completeness we recall an example given by Kuratowski and Sierpiński in \cite{ku-sie}, of a function satisfying the assertion of Proposition \ref{covering}.

\begin{example}\label{kur-sie}
	Let $\psi:[0,1]\to [0,1]$ be given by the formula
	$$
	\psi(t)=\sum_{n=1}^{\infty} \frac{\phi(t-q_n)}{2^n}, 
	$$
	where $\phi(t)=|\sin(\frac{1}{t})|$ for $t\neq 0$ and $\phi(0)=0$, and $(q_1, q_2,\ldots)$ is an injective enumeration of $\QQ\cap [0,1]$.
	
	Then $f=\psi|([0,1]\setminus\QQ)$ satisfies the  assertion of Proposition \ref{covering}. To see this, let us fix $q_n$, and let 
		$$
	\sigma(t)=\sum_{m\neq n} \frac{\phi(t-q_m)}{2^m} \quad\hbox{for}\ t\in [0,1].
	$$
	
	Then $\sigma$ is continuous at $q_n$,
		$$
	\psi(t)=\sigma (t) + \frac{\phi(t-q_n)}{2^n}, 
	$$
	and the definition of $\phi$ yields that
	$$
	\overline{G(f)}\cap (\{q_n\}\times [0,1])=	\{q_n\}\times [\sigma(t), \sigma(t)+2^{-n}].
	 $$ 
\end{example}

\smallskip 

By combining Theorem \ref{d=c implies C_8} with  Proposition \ref{C_8->C_9} and Corollary \ref{b,d}, we immediately get the following corollary. 

\begin{cor}\label{C_8_iff_C_9}
	If the cardinal $\co$ is regular, then the following statements are equivalent:
	\begin{enumerate}
		\item $C_8(\co,\co)$,
		\item $C_9(\co,\co)$,
		\item $\mathfrak{d}=\co$.
		
	\end{enumerate}
\end{cor}

Without any additional set-theoretic assumptions we have  (in ZFC) the following weaker result.

 \begin{theorem}\label{C_9(c,kappa)_Ra_C_8(c,c)}
	For any uncountable cardinal $\kappa<\co$, 
	$$
	C_9(\co,\kappa)\Ra C_8(\co,\co).
	$$ 
\end{theorem}	

\smallskip 

Before giving the proof of Theorem \ref{C_9(c,kappa)_Ra_C_8(c,c)}
let us recall that a subset $S$ of a Polish space $X$ is {\sl 
$\kappa$-concentrated
 in $X$ on a set $D\sub X$}, 
  $\aleph_1\leq \kappa\leq\co$,
  if $|S\setminus U|<\kappa$
  for every open in $X$ set $U\sub X$ that contains $D$. 
  
  \ss 
  
   Let us assume that $C\sub\RR$ is the Cantor set, $Q$ is a countable dense set in $C$ and   $P=C\setminus Q$. 
  \ss 
  
  Let us note that if 
  $S\sub P$, then $S$ is $\kappa$-concentrated on $Q$ in $C$ if and only if $S$ is a $\kappa$-$K$-Lusin in $P$.
 (cf. \cite[Proposition 3.4]{ba-ha}).
 
\smallskip
 

 We shall need the following lemma.
 
 \ss 
 
 \begin{lemma}\label{special H}
 Suppose that there is a set $S\sub P$ of cardinality $\co$ $\kappa$-concentrated on $Q$ in $C$. Then there is a set $H\sub S\times C$  such that for every compact set $K$ in $C\times C$ with $K\cap (Q\times C)$ countable, we have $|H\cap K|<\co$.   
 \end{lemma}
 	
 	\begin{proof}[Proof of Lemma \ref{special H}]
 	We list all compact sets in $C\times C$ with $K\cap (Q\times C)$ countable as $( K_\alpha: \alpha<\co)$.
 	Then we inductively pick distinct points $s_\alpha\in S$ and $t_\alpha\in C$  such that  $(s_\alpha,t_\alpha)\not\in\bigcup\limits_{\beta<\alpha}K_\beta$ for $\alpha<\co$.
 
 	At stage $\alpha<\co$, let us notice that   $\bigcup\limits_{\beta<\alpha}K_\beta \cap (Q\times C)$ has cardinality less than $\co$, and therefore, there is $t\in C$ such that $\bigcup\limits_{\beta<\alpha}K_\beta\cap (Q\times \{t\})=\emptyset$. Since $S\times \{t\}$ is $\kappa$-concentrated on $Q\times \{t\}$ in $C\times\{t\}$, it follows that, for each $\beta<\alpha$, $|K_\beta \cap (S\times \{t\})|<\kappa$ . Therefore, 
 	$$
 	|\bigcup\limits_{\beta<\alpha}\bigl(K_\beta \cap (S\times \{t\})\bigr)|<\co,
 	$$
 	so we can pick $s\in S$, distinct from $s_\beta$ for all $\beta<\alpha$, such that $(s,t)\notin\bigcup\limits_{\beta<\alpha}K_\beta$, and we let $(s_\alpha,t_\alpha)=(s,t)$.
 	
 	Finally, we let $H=\{(s_\alpha,t_\alpha):\alpha<\co\}$.
\end{proof}

With Lemma  \ref{special H} in hand, we proceed with the proof of Theorem \ref{C_9(c,kappa)_Ra_C_8(c,c)} as follows.

\begin{proof}[Proof of Theorem \ref{C_9(c,kappa)_Ra_C_8(c,c)}]
	Let us assume that statement $C_9(\co,\kappa)$ is true. With the help of Theorem \ref{C_9} and the fact that $P$ is a homeomorphic copy of $\baire$, we get a  $\kappa$-$K$-Lusin $S$ in $P$ of cardinality $\co$.   
	
	There is a continuous map $\phi:C\times C\w C$ such that
	$\phi|(P\times C)$ is a homeomorphism  onto $G=\phi(P\times C)$, and the set $D=\phi(Q\times C)$ is countable and disjoint from $G$ (a simple argument to this effect is given in \cite[Lemma 4.2]{mazurkiewicz}). 
	
	 Since $S$ is $\kappa$-concentrated on on $Q$ in $C$ there is, by Lemma \ref{special H}, a set $H\sub S\times C$  such that for every compact set $K$ in $C\times C$ with $K\cap (Q\times C)$ countable, we have $|H\cap K|<\co$. Let $E=\phi(H)$ and $f=\phi^{-1}|E: E\to H$. Upon an embedding of $C\times C$ in $\RR$, we can consider $f$ as a function from a subset of $\RR$ of cardinality $\co$ into $\RR$ and we are going to prove that it is a witness that statement $C_8(\co,\co)$ is true.
	
	Aiming at a contradiction, assume that $f|A$ is uniformly continuous (with respect to any metric compatible with the topology of $C\times C$) on a set $A\sub E$ of cardinality $\co$ and let $B=f(A)=\phi^{-1}(A)$. Then, since $\phi|B:B\to A$ is also uniformly continuous, the function $f|A$ extends to a homeomorphism $\tilde{f}: \bar{A}\to \bar{B}$, where  $\bar{A}$ and  $\bar{B}$ are the closures of $A$ and $B$ in $C$ and $C\times C$, respectively (cf. \cite[Theorem 4.3.17]{eng}). Since $D$ is countable and 
	$\tilde{f}^{-1}=\phi|\bar{B}$ injectively maps $\bar{B}\cap (Q\times C)$ into $\bar{A}\cap D$,
	 $\bar{B}\cap (Q\times C)$ is also countable, and we have $\bar{B}=K_\alpha$ for some $\alpha<\co$. This, however, is impossible, as by Lemma \ref{special H}, we have $|H\cap K_\alpha|<\co$, but on the other hand $B\sub H\cap K_\alpha$ has cardinality $\co$.  
	
\end{proof}

  \smallskip 
  
 A similar reasoning shows, in particular, that also the statements $C_8(\aleph_1,\aleph_1)$ and $C_9(\aleph_1,\aleph_1)$ are equivalent.
 
 \smallskip

Let us  recall that $E$ is a {\sl $\lambda'$-set} in $X$ if for every countable set $L$ in $X$, $L$ is a relative $G_\delta$-set in $E\cup L$; by a theorem of Sierpiński, there is (in ZFC) an uncountable  $\lambda'$-set in $\RR$, cf. \cite{kur}.

\begin{theorem}\label{C_8(w_1,w_1)_iff_C_9(w_1,w_1)} For any uncountable cardinal  $\nu\leq\co$, the existence of a $\lambda'$-set  $T$ of cardinality $\nu$ in the Cantor set $C\sub\RR$ and a $K$-Lusin set $S$ in $P=C\setminus Q$ of cardinality $\nu$, where $Q$ is a countable dense set in $C$, implies $C_8(\nu,\aleph_1)$. Consequently, the existence of a $\lambda'$-set  of cardinality $\nu$ in the Cantor set $C$ implies that  $C_8(\nu,\aleph_1)\Lra C_9(\nu,\aleph_1)$ and hence,
	 we have (in ZFC)
 	$
 	C_8(\aleph_1,\aleph_1)\Lra C_9(\aleph_1,\aleph_1),
 	$ 
 	and each of the statements is equivalent to the assertion $\mathfrak{b}=\aleph_1$.
\end{theorem}
\begin{proof}
Let $H$ be the graph of a bijection from $S$ onto $T$. Since $S$ is concentrated in $C$ on $Q$, it follows that $H$ is  concentrated in $C\times C$ on $Q\times C$. Indeed, if $U$ is an arbitrary open set in $C\times C$ containing $Q\times C$, and $D=(C\times C) \setminus U$, then $V=C\setminus 	\hbox{proj}_1(D)$ is an open set in $C$ containing $Q$  (where $\hbox{proj}_1$ is the projection of $C\times C$ onto the first axis).  Thus $V$ contains all but countably many points of $S$, and consequently, the set $H\setminus U$ is countable.

As in the proof of Theorem \ref{C_9(c,kappa)_Ra_C_8(c,c)}, let $\phi:C\times C\w C$ be a continuous map  such that
 $\phi|(P\times C)$ is a homeomorphism  onto $G=\phi(P\times C)$, and the set $D=\phi(Q\times C)$ is countable and disjoint from $G$.
 
 Let $E=\phi(H)$ and $f=\phi^{-1}|E: E\to H$.
  Upon an embedding of $C\times C$ in $\RR$, we consider $f$ as a function from a subset of $\RR$ of cardinality $\nu$ into $\RR$
  and we shall prove that it is a witness that statement $C_8(\nu,\aleph_1)$ is true.
 
  Aiming at a contradiction, assume that $f|A$ is uniformly continuous on an uncountable set $A\sub E$, let $B=f(A)=\phi^{-1}(A)$ and 
   extend the function $f|A$ to a homeomorphism $\tilde{f}: \bar{A}\to \bar{B}$, where  $\bar{A}$ and  $\bar{B}$ are the closures of $A$ and $B$ in $C$ and $C\times C$, respectively (cf. the proof of Theorem \ref{C_9(c,kappa)_Ra_C_8(c,c)}).

 Let us notice that $E$ is concentrated on $D$ in $C$ and $\bar{A}\cap D\neq\emptyset$. It is easy to see that this implies that $A$ is concentrated on $\bar{A}\cap D$ in $\bar{A}$, hence also $B$ is concentrated on $L=\tilde{f}(\bar{A}\cap D)$ in $\bar{B}$. Clearly, $L$ is a countable subset of $\bar{B}\sub C\times C$ and $B$ is concentrated on $L$ also in $C\times C$. Therefore, $B'= \hbox{proj}_2(B)$ is concentrated on $L'=\hbox{proj}_2(L)$ in $C$ (where $\hbox{proj}_2$ is the projection of $C\times C$ onto the second axis). Since  $H$ is the graph of an injection, $B'$ is uncountable, and  it follows that $L'$ is a countable set in $C$ which is not a $G_\delta$-set in $B'\cup L'$. This, however, contradicts the fact that $T$ is a $\lambda'$-set in $C$ and $B'\sub T$.
 
 \smallskip 
 
 The assertion, stating the equivalence of  statements $C_8(\nu,\aleph_1)$ and $C_9(\nu,\aleph_1)$ assuming the existence of a $\lambda'$-set  of cardinality $\nu$ in the Cantor set follows now from Theorem \ref{C_9}  and Proposition \ref{C_8->C_9}. Indeed, $C_9(\nu,\aleph_1)$ implies that 
 there is also a $K$-Lusin set in $P=C\setminus Q$ of cardinality $\nu$, where $Q$ is a countable dense set in $C$ (cf.  Theorem \ref{C_9}), which by what we have already proved, yields $C_8(\nu,\aleph_1)$. The converse implication is always true (see Proposition \ref{C_8->C_9}).
 
 \smallskip 
 
 The final assertions follow  now from Corollary \ref{b,d}(1) and the existence (in ZFC) of a $\lambda'$-set of cardinality $\aleph_1$ in the Cantor set.
 \end{proof}


While the status of the implication $C_9 \Ra C_8$, the central topic of this note, remains unclear, the following conditions, sufficient for the validity of $C_9 \Ra C_8$, hint at difficulties in finding a model of ZFC where, possibly, $C_9$ is true but $C_8$ is false. 

\begin{proposition}
	If either there 
	are no $K$-Lusin sets in $\baire$ of cardinality $\co$ (in particular, if either $\mathfrak{b}>\aleph_1$ or $\mathfrak{d}<\co$) or at least one of the following statements is true:
	\begin{enumerate}
		\item there exists a Lusin set in $\RR$ of cardinality $\co$,
		\item there exists a $\lambda'$-set  in the Cantor set of cardinality $\co$.
	\end{enumerate}
then $C_8 \Lra C_9$.
\end{proposition}

\begin{proof}
The non-existence of $K$-Lusin sets in $\baire$ of cardinality $\co$ makes $C_9$ false by Theorem \ref{C_9}. Similarly, the existence of a Lusin set in $\RR$ of cardinality $\co$ makes $C_8$ true
(see \cite[proof of Th\'eor\`eme 6 on page 45]{Si}). 

The existence of a $\lambda'$-set  of cardinality $\co$ in the Cantor set $C$ implies that  $C_8(\co,\aleph_1)\Lra C_9(\co,\aleph_1)$, by Theorem \ref{C_8(w_1,w_1)_iff_C_9(w_1,w_1)}.

\end{proof}

\section{Comments}\label{sec:4}
In this section we gathered some additional results related to the subject of this note; their proofs will be presented elsewhere.

\subsection{Mappings into compact spaces}\label{sec:4.1}
One can show that 
 statement  $C_8(\lambda,\kappa)$ is equivalent
to its analogue concerning  functions with ranges in  compact metric spaces.

\begin{proposition}\label{maps_into_compact} 
	For any uncountable cardinals $\kappa\leq\lambda\leq\co$ if the cofinality of $\lambda$ is uncountable, then	the following are equivalent:
	\begin{enumerate}
		\item $C_8(\lambda,\kappa)$,

	\item there is  a set $E\sub [0,1]$ of cardinality $\lambda$, such that for every  uncountable compact metric space $Y$ there is a continuous function $f:E\to Y$, which is not uniformly continuous on any subset of $E$ of cardinality $\kappa$.
	\end{enumerate} 
\end{proposition}
Moreover, if  a set $E$ satisfies property (2) above for a continuous function $f:E\to Y$ into an uncountable compact metric space $Y$, then so does it for a certain real-valued continuous function on $E$. 

However, the latter is no longer true when we replace $E\sub [0,1]$ by an arbitrary separable metrizable space, as shown by the following result.  

\begin{proposition}
	Assuming that no family of less than $\co$ meager sets covers $\RR$, there exists a set of positive dimension $E\sub [0,1]^\NN$ such that 

	\begin{enumerate}
		\item there is a continuous function $f:E\w [0,1]^\NN$, which is not uniformly continuous on any subset of $E$ of cardinality $\co$,
		
		\item each continuous function $g:E\to \RR$ is constant on a subset of $E$ of cardinality $\co$.
	\end{enumerate}
	\end{proposition}

\subsection{Mappings into non-compact spaces}\label{sec:4.2}
One can show that also statement $C_9(\lambda,\kappa)$ is equivalent to an analogue of statement $C_8(\lambda,\kappa)$ but concerning functions with ranges in non-$\sigma$-compact Polish spaces.
	\begin{proposition}\label{maps_to_noncompact}
		For any uncountable cardinals $\kappa\leq\lambda\leq\co$	the following are equivalent:
	\begin{enumerate}
		\item $C_9(\lambda,\kappa)$,
		
			\item	 there is a set $E\sub [0,1]$ 
	of cardinality $\lambda$ such that for every non-$\sigma$-compact Polish space $Y$ there is a continuous function on $E$ which is not uniformly continuous (with respect to any complete metric on $Y$) on any subset of $E$ of cardinality $\kappa$.
\end{enumerate}
\end{proposition}

Moreover, any $\kappa$-K-Lusin set $E$ of cardinality $\lambda$ in $\PP$ has property (2) above.
 In particular, with the help of Theorem \ref{K-Lusin} and  Proposition  \ref{maps_into_compact},  it follows 
that, under 
CH, there exists a $K$-Lusin set $E$ in $\PP$ of cardinality $\co$ such that $E$ admits a continuous function $f:E\to \RR^\NN$ which is not uniformly continuous on any uncountable subset of $E$, but  each continuous map $g:E\to [0,1]^\NN$ is uniformly continuous  on an  uncountable subset of $E$.

\subsection{A characterization of completeness}\label{sec:4.3}
One can show that the existence of a function sequence described in Theorem \ref{seq_on_Polish} characterizes 
completeness of a separable metrizable space $X$.

In fact, the following more general result can be obtained (for terminology see \cite{eng}).

\begin{proposition}
	Let $X$ be a Hausdorff space. Then $X$ is a \v{C}ech-complete Lindel\"{o}f space if and only if there is a sequence $f_1\geq f_2\geq\ldots$ of continuous functions $f_n:X\to [0,1]$ converging pointwise to zero but not converging uniformly on any closed non-compact set in $X$.
\end{proposition}

\bibliographystyle{amsplain}

\end{document}